\documentclass[11pt]{article}
\usepackage{amsmath,amsthm}
\usepackage{amssymb,latexsym}

\def\ba#1{\begin{array}{#1}}
\def\ea{\end{array}}
\def\beq#1{\begin{equation}\label{#1}}
\def\eeq{\end{equation}}
\renewcommand{\Re}{{\rm Re}\,}
\renewcommand{\Im}{{\rm Im}\,}
\newcommand{\RR}{\mathbb{R}}
\newcommand{\Rp}{\RR_+}

\newcommand{\CC}{\mathbb{C}}
\newcommand{\HH}{\mathbb{H}} % upper half plane
\newcommand{\Cp}{\CC_+} % right half-plane
\newcommand{\Cpc}{\overline{\CC}_+} % closed right half-plane

\newcommand{\eps}{\varepsilon}

\newcommand{\Lt}{\mathcal{L}} % Laplace transform
\newcommand{\Pt}{\mathcal{P}} % Poisson integral
\newcommand{\HYp}{\mathrm{HY}(p)} %Hausdorff-Young property of order p
\newcommand{\HYN}[2]{N_{\mathcal{HY},{#1}}(#2)} %Hausdorff-Young norm \HYN{p}{\mu}
\newcommand{\CN}[1]{N_{\mathcal{C}}(#1)} %Carleson norm
\def\dst{\displaystyle}
\def\eop{\hfill$\Box$}
\def\rem#1{}
  
\newtheorem{thm}{Theorem}
\newtheorem{defn}{Definition}

\begin{document}

%---  Frontmatter ---

%%%%% To ease editing, for IMPAN journals add:
\baselineskip=17pt

\title{Characterization of Carleson measures by the Hausdorff-Young property}
\author{Sergey Sadov}
\date{}
\maketitle

%% Classification and key words; note that the 2010 classification is used:
\renewcommand{\thefootnote}{}
\footnote{E-mail: sergey@mun.ca}

\renewcommand{\thefootnote}{\arabic{footnote}}
\setcounter{footnote}{0}

\begin{abstract}
It is shown that the Laplace transform of an $L^p$ ($1<p\leq 2$) function defined on the positive semiaxis
satisfies the Hausdorff-Young type inequality with a positive weight in the right complex half-plane if and only if
the weight is a Carleson measure.
\end{abstract}

%--- End Frontmatter ---

\section{Main theorem}

The classical Hausdorff-Young inequality for the (one-dimensional) Fourier transform 
$$
u(t)\to \mathcal{F}u(\xi)=\int_\RR u(t) e^{-it\xi}\,dt
$$
says
\beq{HY}
 \|\mathcal{F} u\|_{L^{p'}(\RR)}\,\leq\, B(p)\, \|u\|_{L^p(\RR)}.
\eeq
Here $1\leq p\leq 2$, $p'=(1-1/p)^{-1}$ is the conjugate exponent, with $p=1$ corresponding to $p'=\infty$. 
Titchmarsh's now-textbook estimate $B(p)\leq (2\pi)^{1/p'}$ follows from the Parceval theorem and the Riesz-Thorin interpolation theorem.
The sharp constant, which will be of some importance here,   
\beq{Bp}
B(p)=(2\pi)^{1/p'}\left(\frac{p^{1/p}}{{p'}^{1/p'}}\right)^{1/2}
\eeq
has been determined by Babenko \cite{Babenko61}
for even integer $p'$ and Beckner \cite{Beckner75} in the general case.   

We will be dealing with functions $u(t)$ defined on the positive half-line $\Rp=(0,+\infty)$, in which case the Fourier transform is analytic
in the upper half-plane $\HH=\{z\,|\,\Im z>0\}$ and belongs to the Hardy class $H^{p'}(\HH)$ whenever $u\in L^p(\Rp)$, $1\leq p\leq 2$.   

In an equivalent setting, which we prefer for the reason of complex-conjugate symmetry, the upper half-plane $\HH$ is replaced by
the right half-plane $\Cp=\{z\,|\,\Re z>0\}$ and the Fourier transform is replaced by the Laplace transform
\beq{LT}
 \Lt u(z)=\int_0^\infty u(t)\, e^{-zt}\,dt.
\eeq
 
Let $H^s(\Cp)$ denote the Hardy class for the right half-plane. For $s\geq 1$, the norm in $H^s(\Cp)$ of a function $v(z)$ is 
$$
 \|v\|_{H^s(\Cp)}=\sup_{x>0}\,\left(\int_{-\infty}^\infty |v(x+iy)|^{s}\,dy\right)^{1/s}.
$$ 
As a consequence of (\ref{HY}), we have
\beq{HYc}
 \|\Lt u\|_{H^{p'}(\Cp)}\leq B(p) \|u\|_{L^p(\Rp)}.
\eeq

Our main theorem asserts equivalence of the two classes of measures: those with Hausdorff-Young property of order $p>1$, and Carleson measures. 

\begin{defn}
Let $\mu$ be a non-negative Borel measure supported
on the closed right half-plane $\Cpc=\{z\;|\;\Re z\geq 0\}$. 
We say that $\mu$ {\em has the Hausdorff-Young property of order $p$}\ or, in short, that $\mu$ is $\HYp$,
and write $\mu\in \HYp$ if there exists a constant $C$ such that
$$
 \|\Lt u(z)\|_{L^{p'}(\Cpc,d\mu)}
\leq C \|u\|_{L^p(\Rp)}
$$
for any $u\in{L^p(\Rp)}$. 
If $C$ is the smallest such constant and $p>1$, we denote $\HYN{p}{\mu}=C^{p'}$, the ``Hausdorff-Young norm of order $p$'' of $\mu$.  
(We leave $\HYN{1}{\mu}$ undefined.)
\end{defn}

The Lebesgue measure along the $y$-axis, $\delta(x)\otimes dy$, is $\HYp$ for all $p\in [1,2]$ according to (\ref{HY}).
The same is true for the Lebesgue measure along the positive $x$-semiaxis, $dx\otimes \delta(y)$ ($x>0$). The corresponding inequality 
\beq{Hineq}
 \|\Lt\|_{L^{p}(\Rp)\to L^{p'}(\Rp)}\,\leq \,\left(\frac{2\pi}{p'}\right)^{1/p'}
\eeq
belongs to Hardy \cite{Hardy33}. 

Note that any positive  Borel measure on $\Cpc$ is $HY(1)$,
due to the trivial pointwise estimate 
$$
 |\Lt u(z)|\leq \|u\|_{L^1(\Rp)}. 
$$

\medskip
Our definition of Carleson measures will be slightly unconventional (cf.~e.g.\ \cite[\S~I.5]{Garnett81}), to include 
a possible nontrivial mass at the boundary. 
Consider a family of squares adjacent to the boundary of $\Cpc$,
\beq{Qah}
Q_{a,h}=\{z\,|\, \Im z\in [a,a+h], \;\Re z\in[0,h]\} \qquad (a\in\RR,\; h>0).
\eeq

\begin{defn}
A nonnegative Borel measure $\mu$ on $\Cpc$ is a {\em Carleson measure}\ and $\CN{\mu}$ its {\em Carleson norm}\  
if 
$$
 \CN{\mu}=\sup_{a,h}\frac{\mu(Q_{a,h})}{h}
$$ 
is finite.
\end{defn} 

\noindent
{\bf Remark.}\@
The measure $\mu$ in Definition~2 can be split as $\mu=\mu_1+\mu_2$, where $\mu_1=\mu|_{x=0}$ is the boundary part of $\mu$
and $\mu_2$ is a Carleson measure in the conventional sense (except that we work in the right half-plane $\Cp$ instead of upper half-plane $\HH$).
The definition implies that $\mu_1$ is absolutely continuous relative to the Lebesgue measure $dy$.
More precisely, the Radon-Nikodym derivative is bounded: 
\beq{RNder}
\frac{d\mu_1}{dy}\leq \CN{\mu},
\eeq
hence $L^q(d\mu_1)\subset L^q(dy)$ for all $q\geq 1$. 

%\medskip
%Our main result states that Definition~1 with $1<p\leq 2$ and Definition~2 describe the same class of measures.

\begin{thm}[\bf {HY-characterization of Carleson measures}]
The following are equivalent:
\\
(a) $\mu$ is a Carleson measure on $\Cpc$;
\\
(b) $\mu\in \HYp$ for some $p\in(1,2]$;
\\
(c) $\mu\in \HYp$ for all $p\in[1,2]$.
\\[0.3ex]
Moreover, for any $p\in(1,2]$ 
\beq{eqnorm}
A_1(p)^{-1}\CN{\mu}\leq \HYN{p}{\mu}\leq A_2(p)\,\CN{\mu},
\eeq
where $A_1(p)\leq 2^{3/2}p'$ and $A_2(p)\leq 160\pi\,\sqrt{e/p'}$.
These estimates are order-sharp as $p'\to\infty$, that is, $A_1(p)\neq o(p')$ and $A_2(p)\neq o({p'}^{-1/2})$. 
\end{thm}

\noindent{\em Proof}. 
We will first show that the qualitative statement (a) $\Rightarrow$ (c) 
is a simple corollary of Carleson's theorem \cite[Th.~II.3.9]{Garnett81}. The implication (b) $\Rightarrow$ (a) and the left inequality in
(\ref{eqnorm}) will be derived similarly to the proof of the converse part of that same theorem. 
The proof of the right inequality in (\ref{eqnorm}) is put over to Section~2.

\medskip
(a) $\Rightarrow$ (c):
The case $p=1$ is trivial, so assume that $1<p\leq 2$.
If $u\in L^p(\Rp)$, then $\Lt u\in H^{p'}(\Cp)$ and (\ref{HYc}) holds. 
Writing $\mu=\mu_1+\mu_2$ as in Remark after Definition~2, we have: $\|\Lt u\|_{L^{p'}(d\mu_1)}<\infty$ by that Remark, and
$\|\Lt u\|_{L^{p'}(d\mu_2)}<\infty$ by Carleson's theorem \cite[Th.~II.3.9, part (a)$\Rightarrow$(b)]{Garnett81}.

% (c) $\Longrightarrow$ (b) trivially.

\medskip
(b) $\Rightarrow$ (a):
%If $p=2$, then $(2\pi)^{-1/2}\Lt$ is a unitary isomorphism $L^2(\Rp)\to H^2(\Cp)$.
%By part (c)$\Rightarrow$(a) of the Carleson theorem cited above, we conclude that $\mu$ is Carleson.
%
%This argument does not work for $1<p<2$, since $\Lt(L^p(\Rp))$ is a proper subset of $H^{p'}(\Cp)$.
%Instead of quoting the assertion of Carleson's theorem, 
%
We mimic the proof of part (c)$\Rightarrow$(a) of the Carleson theorem cited above.
Let $u\in L^p(\Rp)$ be a function with $\|u\|_p=1$ whose Laplace transform $v=\Lt u$ satisfies $b=\inf_{z\in Q_{0,1}} |v(z)|>0$.      
For $h>0$, define
$u_h(t) = h^{1/p} u(ht)$. 
Then $\|u_h\|_p=1$ and 
$
 \;v_h(z)=\Lt u_h(z)=h^{-1/p'} v(z/h).
$
If $z\in Q_{0,h}$, then  
$\,
 |v_h(z)|\geq h^{-1/p'} b
$.
Therefore 
$$
 \|v_h\|^{p'}_{L^{p'}(d\mu)}\geq \int_{Q_h} |v_h(z)|^{p'}\,d\mu(z) \geq \mu(Q_h) h^{-1} b^{p'}.
$$ 
On the other hand, by condition (b) of the Theorem,
$$
 \|v_h\|^{p'}_{L^{p'}(d\mu)}\leq \HYN{p}{\mu} \,\|u_h\|^{p'}_{p} =\HYN{p}{\mu}.
$$ 
It follows that
$$
 \mu(Q_{0,h})\leq  h b^{-p'} \HYN{p}{\mu}.
$$
Estimates for $\mu(Q_{a,h})$ for $a\neq 0$ are obtained similarly by considering the test functions $u_{h}(t) e^{iat}$. 
We conclude that $\mu$ is Carleson, with Carleson norm estimate 
$\CN{\mu}\leq b^{-p'} \HYN{p}{\mu}$. The left part of the inequality (\ref{eqnorm}) follows with $A_1(p)\leq b^{-p'}$.   

The final part of the argument is aimed at obtaining an estimate for $b^{-p'}$ that grows linearly with $p'$.
Consider the test function 
$$
u(t)=\left\{\ba{l}\eps^{-1/p}, \quad 0<t<\eps,\\ 0,\quad\;\;t>\eps\ea, \right.
$$ 
Then $v(z)=\eps^{1/p'}(1-e^{-z\eps})/(z\eps)$.

The inequality $|e^{-w}-w+1|\leq |w|^2/2$ is valid whenever $\Re w\geq 0$ (say, by Taylor's formula with remainder in the integral form). 
Thus $|(e^{-w}-1)/w|\geq 1-|w|/2$.
Consequently,  
$$
 |v(z)|^{p'}\geq  \eps \left(1-\frac{\eps |z|}{2}\right)^{p'}.
$$
For $z\in Q_{0,1}$, the minimum occurs at the corner $z=1+i$ where $|z|=\sqrt{2}$.  
Finally, setting $\eps=\sqrt{2}/p'$, we get  
$$
b^{-p'}\leq \frac{p'}{\sqrt{2}}\left(1-\frac{1}{p'}\right)^{-p'}\leq \frac{4p'}{\sqrt{2}},
$$ 
as claimed.

\section{Evaluation of constants in the inequalities}

To prove the upper bound for the constant $A_2(p)$ in the right inequality (\ref{eqnorm}),
we need Carleson's estimate for the Poisson integral with a numeric constant.
An independent derivation of it is given below.
(We make no claim that the obtained constant is better than what can reconstructed from proofs in existing literature.)
In the %``automatic'' 
interpolatory part of Carleson's theorem, we use a combination of Marcinkiewicz and Riesz-Thorin theorems to ensure the 
desired behaviour of the constant as a function of $p$.
  
In this theorem, the conventional Carleson measures in $\HH$ (no mass at the boundary) are used. The squares
$Q_{a,h}$ in Definition~2 are to be substituted by the squares $R_{a,h}=\{x+iy\,|\, a< x< a+h,\; 0<y<h\}$ and the formula for
$\CN{\mu}$ is to be modified accordingly.  

In addition, the following notation will be needed:
\\
$\diamond\;\;\;$
$P_a(t)=\pi^{-1}a/(a^2+t^2)$, the Poisson kernel;
\\
$\diamond\;\;\;$
$\Pt$, the Poisson convolution operator for the upper half-plane: 
$$
 \Pt f(x+iy)=\int_{-\infty}^\infty P_y(x-t)f(t)\,dt;
$$ 
$\diamond\;\;\;$
$E_g(\lambda)$, the ``large value set'': for a function $g(z)$ and $\lambda>0$, 
$$
E_g(\lambda)=\{z\,:\, |g(z)|>\lambda\}.
$$ 

\begin{thm}[\bf {Carleson's theorem with numeric constant}] 
Let $\mu$ be a Carleson measure in $\HH$.
Suppose that $f(x)\in L^p(\RR)$ and $g(x+iy)=\Pt f(x+iy)$.
\\[0.3ex]
(a)  If $p=1$, then 
\beq{Car_1}
 \mu(E_g(\lambda))\leq 10\, \CN{\mu}\,\|f\|_1.
\eeq
(b) If $1<p<\infty$, then
\beq{Car_p}
 \|g\|_{L^p(d\mu)} \leq \left(M(p)\,\CN{\mu}\right)^{1/p}\,\|f\|_p,
\eeq
where $M(p)\leq 40 p'$ when $1<p<2$, and $M(p)\leq 79$ when $p\geq 2$.
\end{thm}

\medskip
We will finish the proof of Theorem~1 and prove Theorem~2 afterwards.
 
Given a Carleson measure $\mu$ on $\Cpc$, write $\mu=\mu_1+\mu_2$ as in the Remark after Definition~2.
For a function $v\in H^{p'}(\Cp)$, we have
$$
\|v\|_{L^{p'}(d\mu)}^{p'}=\|v\|_{L^{p'}(d\mu_1)}^{p'}+
\|v\|_{L^{p'}(d\mu_2)}^{p'}\leq \CN{\mu}\,(1+ M(p')) \|v\|_{H^{p'}}.
$$ 
by (\ref{RNder}) and (\ref{Car_p}).
It follows by (\ref{HYc}) that for $v=\Lt u$ 
$$
 \|v\|_{L^{p'}(d\mu)}^{p'}\leq (1+M(p'))\cdot (B(p))^{p'}.
$$
(The appearance of $M+1$ explains why we have favored the mingy constant 79 over a generous 80 in Theorem~2.)\@ 
Substituting the evaluation of $B(p)$ from (\ref{Bp}) and the estimate for $M(p')$ from Theorem~2, we 
get
$$
\HYN{p}{\mu}\leq 80\cdot 2\pi\, p^{p'/2p}\, p'^{-1/2}.
$$
To obtain the upper bound for $A_2(p)$ as claimed in Theorem~1, it remains to notice that $\sup_{1<p\leq 2}(p'p^{-1}\ln p) =\lim_{p\to 1}(\dots)=1$.
   
\smallskip
Finally, let us show that the obtained upper bounds for $A_1(p)$ and $A_2(p)$ in (\ref{eqnorm}) are order-sharp.

\smallskip
1. Take $\mu=\delta(z-1/p')$, that is, $\int v(z)\,d\mu=v(1/p')$.
Clearly, $\CN{\mu}=p'$, while H\"older's inequality implies that $\HYN{p}{\mu}\leq 1$. Hence $A_1(p)\geq p'$.

\smallskip
2. Take $\mu=dy$, the Lebesgue measure along the imaginary
axis. Then $\CN(\mu)=1$, while $\HYN{p}{\mu}=B(p)^{p'}$. Hence $A_2(p)\geq B(p)^{p'}\sim cp'{-1/2}$ as $p'\to\infty$.

\smallskip
The proof of Theorem~1 is complete.

\bigskip
\noindent
{\it Proof of Theorem~2}. 
Let us first derive part (b) from part (a).
The inequality (\ref{Car_1}), the trivial inequality $\|g\|_\infty\leq \|f\|_\infty$, and the Marcinkiewicz interpolation theorem 
\cite[Th.~I.4.5]{Garnett81} % esp. inequality on the bottom of p.~27
% or \cite[Th.~1.3.2, (1.3.11)]{Grafakos}. 
yield
$$
 \|g\|_{L^p(d\mu)} \leq 2\left(10 p'\,\CN{\mu}\right)^{1/p}\,\|f\|_p.
$$
We wish to obtain a constant that behaves like $O(1)^{1/p}$ as $p\to\infty$, while the constant in the above inequality
tends to the limit $2$. 
To optimize the upper bound, note that by the Riesz-Thorin theorem
$$
 \left(\frac{\|g\|_{L^p(d\mu)}}{\|f\|_p}\right)^p\; \leq\; 
\inf_{1<r\leq p}  \left(\frac{\|g\|_{L^r(d\mu)}}{\|f\|_r}\right)^r\;\leq\;
\inf_{1<r\leq p} \left(2^r \cdot 10 r'\,\CN{\mu}\right).
$$
The function $r\,\to\,2^r\, r'= 2^r\,r/(r-1)$ attains its minimum $m$ at the root $r_0$
of the equation $r_0(r_0-1)\ln 2=1$. Calculation gives $r_0\approx 1.80104$ and $m\approx 7.83495$.
The upper bounds for $M(p)$ as stated in Theorem~2 are obtained by rounding up.    
   
\medskip
Our proof of part (a) is a shortcut of a standard proof \cite{Garnett81}. 
The three underlying steps are:
Calderon-Zygmund decomposition 
$\Rightarrow$ Hardy-Littlewood maximal theorem 
$\Rightarrow$ Estimate for nontangent maximal function 
$\Rightarrow$ Carleson's theorem. 
These steps will be implicit in our calculation.

Fix $\lambda>0$ and consider the Calder\'on-Zygmund decomposition for $f$ at height $\alpha=\lambda/7$ \cite[Lemma~VI.2.2]{Garnett81}: 
$\RR=B\cup G$, $B\cap G=\emptyset$, 
$f(x)\leq \alpha$ for almost every $x\in G$, $B=\bigcup I_j$, a finite or countable union of disjoint intervals, and 
$$
\alpha\leq M_j \leq 2\alpha,
$$ 
where
$$
M_j=\frac{1}{|I_j|}\int_{I_j} |f| dx.
$$
We may assume for simplicity that the number of intervals $I_j$ is finite: functions $f$ for which this is true are dense in $L^1$.

If $I_j=[a_j,b_j]$, let $\tilde I_j=\big[a_j-|I_j|,\,b_j+|I_j|\big]$, so that $|\tilde I_j|=3|I_j|$. 
Let $R_j$ be the square in $\HH$ with base $\tilde I_j$, i.e.\
$x+iy\in R_j\;\Leftrightarrow x\in \tilde I_j,\; 0<y\leq 3|I_j|$.
Our goal is to show that $E_g(\lambda)\subset \cup R_j$.

\smallskip
Fix $\,z=x_0+iy\notin \cup R_j\,$ and consider separately contributions to $g(x_0+iy)$ of $f$ restricted to the intervals $I_j$
according to whether $\tilde I_j$ contains $x_0$ or lies to the right, resp.\ to the left of $x_0$. Formally: let $x_j$ be the point
in $I_j$ closest to $x_0$. 
Define the mutually disjoint sets of indices $S_0$, $S_+$, $S_-$ as follows: $j\in S_0$ (resp.\ $S_+$ or $S_-$) if $x_j-x_0=0$ (resp.\ $>0$ or $<0$). 
%$$
%\ba{l}
%S_0=\{j\,|\, x_0\in \tilde I_j\}, \\
%S_+ =\{j\,|\, x_0< \inf\tilde I_j\},\\
%S_- =\{j\,|\, x_0> \sup\tilde I_j\}.
%\ea
%$$

In case $S_0\neq\emptyset$, let $L=\max_{j\in S_0}|I_{j}|$ occur for $j=j_0$.   
Then $y>3L$, as otherwise we would have $z\in R_{j_0}$. Clearly, 
$I_j\subset [x_0-2L, x_0+2L]$ for any $j\in S_0$, hence $\sum_{j\in S_0} |I_j|\leq 4L$.
Since $P_y(t)<(\pi y)^{-1}$, we get
\beq{PS0}
 \sum_{j\in S_0}\int P_{y}(x_0-t) |f(t)|\,dt \leq \frac{1}{\pi y} \sum_{j\in S_0} M_j|I_j|\leq \frac{2\alpha\cdot 4L}{\pi y}
<\frac{8}{3\pi}\alpha.
\eeq
 
Let us now evaluate contribution of the intervals $I_j$ with $j\in S_+$.
Define the counting function for the total length of such intervals: 
$$
 F(t)=\sum |I_j| \quad \mbox{\rm over $j$ such that $x_0+|I_j|<a_j\leq x_0+t$}. 
$$
The function $F$ defined in $(0,+\infty)$ is nondecreasing, upper-semicontionuos, and $F(x)=0$ in the right neighborhood of $x_0$.
In addition, we have the important inequalities
\beq{F2}
 F(b_j-x_0)\leq b_j-x_0 \quad (j\in S_+) 
\qquad\mbox{\rm and}\quad
 F(x)< 2x. 
\eeq
The first inequality is obvious; moreover, if $b_{j}\leq x<a_{j+1}$, then $F(x-x_0)<x-x_0$. 
And if $a_{j}\leq x<b_{j}$, then $F(x-x_0)\leq (a_j-x_0)+|I_j|< 2(a_j-x_0)$. 
%, since the condition $x\notin \tilde I_j$ implies $a_j-x_0>|I_j|$.

Let $J=\max j$. We have
$$
\ba{l}
\dst
\sum_{j\in S_+} \int_{I_j} P_y(x_0-t) |f(t)|\,dt
\\[4ex]\dst
\leq
\sum_{j\in S_+} P_y(a_j-x_0) \int_{I_j}|f(t)|\,dt
\\[4ex]\dst
\leq
2\alpha \sum_{j\in S_+} P_y(a_j-x_0) |I_j|
\\[4ex]\dst
=
2\alpha \int_{0}^{b_J-x_0} P_y(t) dF(t).  
\ea
$$
Integrating by parts, using (\ref{F2}) and the elementary inequality $tP_y(t)\leq (2\pi)^{-1}$, we get
$$
\ba{l}
\dst
\sum_{j\in S_+} \int_{I_j} P_y(x-t) |f(t)|\,dt
\\[4ex]\dst
\leq
2\alpha \left(\frac{1}{2\pi}-\int_{0}^{\infty} \frac{d}{dt}P_y(t) \,2t\,dt\right)  
\\[4ex]\dst
\leq
2\alpha\left(\frac{1}{2\pi}+1\right). 
\ea
$$
The contribution of the intervals $I_j$ with $j\in S_-$ has the same upper bound.
Combining with (\ref{PS0}), we obtain
$$
\int_B P_y(x_0-t) |f(t)|\,dt \leq  \alpha\left(\frac{2}{\pi}+4+\frac{8}{3\pi}\right) < \frac{11}{2}\alpha.
$$
Finally,
$\sup_G |f|\leq \alpha$, and we conclude:
$$
|g(x_0+iy)|\leq \int_{B\cup G}\, P_y(x_0-t) |f(t)|\,dt  \, < \frac{13}{2}\alpha.
$$ 

\smallskip
Summarizing, we can cover the set $E_g(\lambda)$ by the union of squares $R_{j}$. The total of their sidelengths
is $\sum 3|I_j|\leq 3\|f\|_1/\alpha<10 \|f\|_1/\lambda $.
The inequality (\ref{Car_1}) follows.
\eop

\subsection*{Acknowledgements}

The research that led to this paper was considerably driven by a collaboration with Professor Anatoli Merzon  
(Universidad Michoacana de San Nicol\'as de Hidalgo, Morelia, M\'exico). The present paper supersedes
\cite{SM-CMA10} and partly supersedes the unpublished joint work \cite{MerzonSadov11}.       

I acknowledge a financial support by NSERC grant during my employment at Memorial University of Newfoundland, Canada,
and a partial travel support from Professor Merzon's CONACYT grant during my visits to Morelia in 2007 and 2009.

%\newpage


\begin{thebibliography}{8}

\bibitem{Babenko61} 
K.\ I.\ Babenko, An inequality in the theory of Fourier integrals, {\em Izv.~Akad.~Nauk~SSSR, Ser.~Mat.}{\bf 25} (1961), 531--542; 
English transl., {\em Amer.~Math.~Soc.~Transl.} (2) {\bf 44}, 115--128.

\bibitem{Beckner75} 
W.~Beckner, {Inequalities in Fourier analysis}, {\em Ann.~of~Math.} {\bf 102} (1975), 159--182.

\bibitem{Garnett81} J.\ B.\ Garnett, {\it Bounded Analytic Functions}, Academic Press, 1981.

\bibitem{Hardy33}
G.\ H.~Hardy, The constants of certain inequalities, 
{\em J.~London Math.~Soc}. {\bf 8} (1933), 114--119.

\bibitem{MerzonSadov11}
A.\ Merzon, S.\ Sadov, Hausdorff-Young type theorems for the Laplace transform restricted to a ray or to a curve in the complex plane, {\tt http://arxiv.org/math.CA/1109.6085}

\bibitem{SM-CMA10} S.\ Sadov, A.\ Merzon,
$L^2$-estimates for the Laplace transform along a family of hyperbolas in the right half-plane,
in: Proceedings of ``Analysis, Mathematical Physics and Applications''
(Ixtapa, Mexico, March 1--5, 2010), 
{\em Comm.\ in Math.\ Analysis},
Conference 03 (2011), 204--208.

\end{thebibliography}
\end{document}